\documentclass[12pt]{article}
\title{Proofs of Two Conjectures by Mecke for Mixed Line-Generated Tessellations}
\author{Eike Biehler, \\Friedrich-Schiller-Universit{\"a}t Jena, \\eikebiehler@web.de}
\date{\today}
\usepackage[latin1]{inputenc}
\usepackage{makeidx}
\usepackage{graphicx}
\usepackage{amssymb}
\usepackage{multirow}
\usepackage[paper=a4paper,left=30mm,right=30mm,top=15mm,bottom=25mm]{geometry}
\usepackage{lscape}
\begin{document}
\newtheorem{Definition}{Definition}
\newtheorem{Satz}{Satz}
\newtheorem{Lemma}{Lemma}
\newtheorem{Korollar}{Korollar}
\newtheorem{Bemerkung}{Bemerkung}
\newtheorem{Vermutung}{Vermutung}

\maketitle

\begin{abstract}
For a compact and convex window, Mecke described a process of tessellations which arise from cell divisions in discrete time. At each time step, one of the existing cells is selected according to an equally-likely law. Independently, a line is thrown onto the window. If the line hits the selected cell the cell is divided. If the line does not hit the selected cell nothing happens in that time step.\\
With a geometric distribution whose parameter depends on the time, Mecke transformed his construction into a continuous-time model. He put forward two conjectures in which he assumed this continuous-time model to have certain properties with respect to their iteration. These conjectures lead to a third conjecture which states the equivalence of the construction of STIT tessellations and Mecke's construction under some homogeneity  conditions.\\
In the present paper, the first two conjectures are proven. A key tool to do that is a property of a continuous-time version of the \textit{equally-likely} model classified by Cowan.
\end{abstract}


\section{Introduction}
The topic of this article are random tessellation processes in the plane. In general, random tessellations are constructed by lines or line segments that are thrown onto the plane under a certain probability law. In our context, line segments are always intersections of lines and a so called cell within a compact and convex window. Both the throwing of the lines and the selection of the cell to be divided are governed by specific probability laws. The timing of the cell division may depend on the selection rule for the cell to be divided.\\
In \cite{Mecke-inhomogen}, Mecke developed a new process in discrete time in a convex and compact window $W$: At the first time step, a line is thrown onto the window according to a law $Q$ dividing the window in two cells almost surely. At the second time step, one of the cells is selected for division according to an equally-likely law. Independently, a line is thrown onto the window. If the line intersects the selected cell, that cell is divided into another two cells. If, however, the line does not intersect the selected cell, although \textit{no} new cell is created, another so-called quasi-cell arises. In any case, the number of quasi-cells (\textit{real} cells plus empty quasi-cells) is always the number of time steps passed plus one.\\
At each time step one of the quasi-cells (which if they are empty cannot actually be divided) is selected with a probability equal to any other cell. If the line thrown independently hits a \textit{real} cell, that cell is divided into two \textit{real} cells. If a \textit{real} cell is selected but the line does not hit it, the \textit{real} cell remains; one new empty quasi-cell is added. If an empty quasi-cell is selected, there automatically arise two new empty quasi-cells.\\
Mecke proposed a way to transfer this process from discrete to continuous time by assigning to an arbitrary time $t$ a geometrically-distributed random number of steps (dependent on $t$) in the discrete process. After formulating two conjectures, he examined the special case of a homogeneous line measure to be used for the (potential) cell division for which he stated another conjecture that his model in continuous time has the same distribution as the STIT tessellation process introduced and examined by Mecke, Nagel and Weiß in e.g. \cite{GlobalConstruction}, \cite{Nagel:2005} and \cite{Mean-values-3D}.\\
While his last conjecture, Conjecture 3 regarding the homogeneous case, was proven in \cite{Conjecture-3} in rather lenghty terms, a by-product of that paper was a way to actually understand Mecke's construction as a \textit{process} in continuous time. By this new access however, which is related to the \textit{equally-likely} model Cowan examined in \cite{Cowan}, the proofs of Mecke's remaining conjectures could be undertaken.\\
In this paper, after a short introduction into Mecke's construction (section \ref{sec: Mecke-process}), the distribution of the lifetime beyond an arbitrary point in time of a convex set within a cell of a fixed tessellation in Mecke's continuous-time model is calculated (section \ref{sec: Lifetime-of-convex-set-in-Mecke-continuous-time}). This allows the proofs in section \ref{sec: Proofs-Conjectures-1-2}.

\section{The Mecke process}
\label{sec: Mecke-process}
Throughout this paper, we will consider, in the Euclidian plane, a compact and convex polygon $W \subset \mathbb{R}^2$ with non-empty interior. Let $[\mathcal{H}, \mathfrak{H}]$ be the measurable space of all lines in $\mathbb{R}^2$ where the $\sigma$-algebra is induced by the Borel $\sigma$-algebra on a parameter space of $\mathcal{H}$. For a set $A \subset \mathbb{R}^2$ we define $$[A] = \{g \in \mathcal{H}: g \cap A \neq \emptyset\}.$$ Let $Q$ be a non-zero locally finite measure on $[\mathcal{H}, \mathfrak{H}]$ which is not concentrated on one direction but which is bundle-free, i.e. there is no point $x \in \mathbb{R}^2$ such that $Q([\{x\}]) = Q(\{g \in \mathcal{H}: g \cap \{x\} \neq \emptyset\}) > 0$. For $Q$, $0 < Q([W]) < \infty$ is true.\\
\subsection{Mecke's process in discrete time}
\label{sec: Mecke-process-discrete}
Let there be lines $\gamma_j, j=1, 2, ...$, that are i.i.d. according to the law $Q([W])^{-1} Q(\cdot \cap [W])$. Further let us use, independently of $\gamma_j$, independent $\alpha_j, j=1, 2, ...$ where $\alpha_j$ is uniformly distributed on the set $\{1, ..., j\}$.\\
If a line $\gamma_j$ does not contain the origin $o$ then $\gamma_j^+$ shall be the open halfplane bounded by $\gamma_j$ which contains the origin. Correspondingly, $\gamma_j^-$ is the open halfplane bounded by $\gamma_j$ which does not contain the origin. As the distribution of $\gamma_j$ is bundle-free, we can neglect the possibility of $\gamma_j$ going through the origin as the probability of this is zero.\\
Let be $\tilde{C}_{0,1}=W$, $\tilde{C}_{1,1}=W \cap \gamma_1^+$ and $\tilde{C}_{1,2}= W \cap \gamma_1^-$. For $n=2, 3, ...$ we define $$\tilde{C}_{n,j} = \left\{\begin{array}{cl}\tilde{C}_{n-1, j} & \textrm{ if $j \in \{1, ..., n\}, j \neq \alpha_n$}\\&\\
\tilde{C}_{n-1, \alpha_n} \cap \gamma_n^- & \textrm{ if $j=\alpha_n$}\\&\\
\tilde{C}_{n-1, \alpha_n} \cap \gamma_n^+ & \textrm{ if $j=n+1$}\end{array}\right.$$
These entities $\tilde{C}_{n, j}$ are called \textbf{quasi-cells}. Some of these quasi-cells are empty. Those quasi-cells that are not empty will be called \textbf{cells}.\\
From this, we can deduce a random process: After each \textbf{decision time} $n$, $n = 1, 2, ...,$ we consider the tessellation $\mathcal{T}_n$ consisting of the quasi-cells $\tilde{C}_{n, 1}, ..., \tilde{C}_{n, n+1}$. This decision time is called the $n$-th decision time accordingly. If, at that decision time, the number of cells (i.e. non-empty quasi-cells) actually changes, that decision time is called a \textbf{jump time}. Obviously, the $k$-th jump time is that decision time at which the number of cells reaches $k+1$. Let us denote the random closed set of the closure of the union of cell boundaries that are not part of the window's boundary at a step $n$ for the tessellation $\mathcal{T}_n$ as $$Y^M_d(n,W) = \overline{\bigcup_{j=1}^{n+1} \partial \tilde{C}_{n, j} \setminus \partial W}.$$\\
Then $(Y^M_d(n,W): n \in \mathbb{N})$ is called the \underline{M}ecke process in \underline{d}iscrete time. Here, $\mathbb{N}=\{0, 1, 2, ...\}$ is the set of the natural numbers.

\subsection{The Mecke model in continuous time}
\label{sec: Mecke-continuous-time}
In \cite[Section 4]{Mecke-inhomogen}, Mecke introduces a mixed line-generated tessellation model such that the tessellation $\mathcal{T}^t$ at the continuous time $t \in [0, \infty)$ corresponds to the tessellation $\mathcal{T}_{\nu(t)}$ at the discrete random time $\nu(t)$ where for the distribution of $\nu(t)$ $$\mathbb{P}(\nu(t) = k) = e^{-t} \left(1-e^{-t}\right)^k, k = 0, 1, ...$$ holds. Mecke used $Q([W])=1$ for his considerations. For general $Q$, i.e. where $Q([W])=1$ is not necessarily true any more, the distribution is
\begin{equation}\label{eq: Allgemeine-Entscheidungszeit-Verteilung}\mathbb{P}(\nu(t) = k) = e^{-Q([W])t} \left(1-e^{-Q([W])t}\right)^k, k = 0, 1, ...\end{equation} This is the geometric distribution with parameter $e^{-Q([W])t}$; the model (which yields a random tessellation for any fixed time $t$ but cannot yet be described as a process) thus has the characteristics $Q$ and $t Q([W])$. (In Mecke's paper, these characteristics were $Q$ and $t$ due to $Q([W])=1$. Here, to have a connection between the characteristic and the exponential function's exponent, the characteristic is called $t Q([W])$.) A possible interpretation is that the decision times are no longer at equidistant discrete times $n=1, 2, ...$ Instead, the law describes how many decisions take place until the time $t$. The $\nu(t)$ are assumed independent of all other random variables that are used in the construction of the Mecke process.

\subsection{The sum of exponentially-distributed random variables}
While there are more general results for the distribution of a sum of exponentially-distributed random variables with unequal parameters (e.g. see \cite{Faltung-Summe-exponentialverteilter-Zufallsgroessen}), for the  special case needed here the following calculations allow a quick understanding. If a random variable $X$ is exponentially distributed with parameter $\lambda$, we will write $X \sim \mathcal{E}(\lambda)$.
\label{sec: Prozess-Eigenschaften}
\begin{Lemma}
\label{Lemma: n-te-Sprungzeit-equally-likely}
Let $n \in \mathbb{N} \setminus \{0\}$ be fixed. Let further $S_n = \sum_{j=1}^n T_j$ be the sum of independent exponentially distributed random variables $T_1, ..., T_n$ with $T_j \sim \mathcal{E}(j R)$ for $j=1, 2, ..., n$ and a fixed $R >0$.
Then $$\mathbb{P}(S_n \leq t) = \int_0^t n R e^{-nxR}(e^{xR}-1)^{n-1}dx = e^{-ntR}(e^{tR}-1)^n = (1-e^{-tR})^n$$ holds.
\end{Lemma}
\textbf{Proof}\\
The proof is by induction. For $n=1$, obviously$$\mathbb{P}(S_1 \leq t) = \mathbb{P}(T_1 \leq t) = \int_0^t R e^{-xR} dx = \int_0^t 1 R \cdot e^{-1 \cdot xR} \cdot (e^{xR}-1)^0 dx =  1-e^{-tR}$$ holds which is true according to the condition $T_1 \sim \mathcal{E}(R)$.\\
Let the lemma be true for $n$. Then, because of $S_{n+1} = S_n + T_{n+1}$ with $T_{n+1} \sim \mathcal{E}\left((n+1)R\right)$ and the independence of $S_n$ and $T_{n+1}$,  for the density of $S_{n+1}$  $$\begin{array}{rl} f_{S_{n+1}}(x) = & f_{S_n+T_{n+1}}(x)\\&\\ = & \int_0^x f_{S_n}(u) f_{T_{n+1}}(x-u) du\\&\\
= & \int_0^x nR e^{-nuR}(e^{uR}-1)^{n-1} (n+1) R e^{-(n+1)(x-u)R}du\\&\\
= & (n+1)Re^{-(n+1)xR} \int_0^x nR e^{uR} (e^{uR}-1)^{n-1}du\\&\\
= & (n+1)Re^{-(n+1)xR} [(e^{uR}-1)^n]_{u=0}^{u=x}\\&\\
= & (n+1)Re^{-(n+1)xR} (e^{xR}-1)^n
\end{array}$$ holds. Integration yields the second equation, straightforward calculation the third equation in the lemma. \hfill $\Box$\\

\begin{Lemma}
\label{Lemma: X-t-in-equally-likely}
Let $N_t = \max\{n: \sum_{j=1}^n T_j \leq t\}$ denote the number of $T_j \sim \mathcal{E}(jR)$, $j=1, 2, ...$, which have consecutively expired until the time $t$. Then for  $k=0,1, 2, ...$
\begin{equation}\label{eq: P-N-t-gleich-k}\mathbb{P}(N_t = k) = e^{-tR} \left(1-e^{-tR}\right)^k\end{equation} holds.
\end{Lemma}
\textbf{Proof}
From the distribution of the $S_k, k=1, 2, ...,$ one gets
$$\begin{array}{rl} & \mathbb{P}(N_t = k)\\&\\ = & \mathbb{P}(S_k \leq t < S_{k+1})\\&\\ = & \mathbb{P}(S_k \leq t) - \mathbb{P}(S_{k+1} \leq t)\\&\\
= & \left(1-e^{-tR}\right)^k - \left(1-e^{-tR}\right)^{k+1}\\&\\
= & \left(1-e^{-tR}\right)^k \left(1 - (1 - e^{-tR})\right)\\&\\
= & e^{-tR} \left(1-e^{-tR}\right)^k
\end{array}$$
For $N_t=0$, the result follows from Lemma \ref{Lemma: n-te-Sprungzeit-equally-likely} immediately.\hfill $\Box$\\

\subsection{The Mecke process in continuous time}
Comparing the equations (\ref{eq: Allgemeine-Entscheidungszeit-Verteilung}) and (\ref{eq: P-N-t-gleich-k}), we see that with $N_t = \nu(t)$ and $R = Q([W])$ both yield the same result. Therefore, the $T_j$ from Lemma \ref{Lemma: X-t-in-equally-likely} with $T_j \sim \mathcal{E}(j Q([W]))$ can be interpreted as the (continuous-time) waiting times for the quasi-state of the tessellation to change from a quasi-state with $j$ quasi-cells to a quasi-state with $j+1$ quasi-cells.
Thus, we can define
\begin{Definition}
\label{Definition: Mecke-continuous-time-process}
Let us have a window $W \subset \mathbb{R}^2$. Let $(Y^M_d(n,W): n \in \mathbb{N})$ be the Mecke process in discrete time as described in section \ref{sec: Mecke-process-discrete}. Let $(N_t: t \geq 0)$ be the process of the number of expired random variables $T_j \sim \mathcal{E}(j Q([W]))$ as in Lemma \ref{Lemma: X-t-in-equally-likely}. Then for every $t \in [0, \infty)$ we define $$Y^M_c(t,W)= Y^M_d(N_t, W)$$
and the \underline{M}ecke process in \underline{c}ontinuous time as $(Y^M_c(t,W): t \geq 0).$\\
\end{Definition}

\section{The waiting time until a convex set is hit in the Mecke process in continuous time}
\label{sec: Lifetime-of-convex-set-in-Mecke-continuous-time}
We will now give a formula for the waiting time of a convex set within a cell in the Mecke process in continuous time to be hit by a line.\\
Let us have a fixed time $s$. We work on the condition that, at this time, the tessellation has $n$ quasi-cells, thus $\mathcal{T}^s = \mathcal{T}_{n-1}$. For the waiting time $T^M_n$ in this state, $T^M_n \sim \mathcal{E}(n Q([W]))$ holds. In this fixed tessellation with $n$ quasi-cells, let us have $\kappa$ ('real') cells.\\
Let these cells be called $C_1, ..., C_\kappa$. Let $S_j \subset C_j \subset W, j=1, ..., \kappa,$ be a convex set within a cell $C_j$, created deterministically from $C_j$ (thus $S_j$ depends on the cell $C_j$!). Let us denote by $X_{S_j}$ the waiting time for such a set $S_j$ to be hit by a line for the first time after the time $s$. Examples for the deterministic function creating $S_j$ from $C_j$ may be the identical function or an intersection of $C_j$ with a fixed convex set. Some of the $S_j$ may be empty which does not compromise the following consideration. Nonetheless, when we consider an $S_j$ in the future it is always assumed to be non-empty.\\
It may be possible that the cell $C_j$ is hit by a line but the set $S_j$ is not. In this case, the waiting time for $S_j$ to be hit shall not begin anew but rather be extended until it is actually hit. The waiting time until the set $S_j$ is hit is the waiting time $T^M_n$ if and only if the cell $C_j$ that contains $S_j$ is selected for division in this step (i.e. $\alpha_n = j$ in Mecke's construction) \textit{and} the set $S_j \subset C_j$ is hit by the line. The probability for this to happen is $$\mathbb{P}(\textrm{$j=\alpha_n$, $S_j \cap \gamma_n \neq \emptyset$}|S_j \subset C_j \in \mathcal{T}_{n-1}) = \frac{1}{n} \frac{Q([S_j])}{Q([W])}.$$\\
If the set is not hit (which happens with probability $1-\frac{1}{n} \frac{Q([S_j])}{Q([W])}$) the waiting time for the set to be hit is the sum of the waiting times $T^M_n + T^M_{n+1}$ if and only if the waiting times $T^M_n$ and $T^M_{n+1}$ have passed and the set is hit in the $(n+1)$-th division step the probability of which is $$\mathbb{P}(\textrm{$j=\alpha_{n+1}$, $S_j \cap \gamma_{n+1} \neq \emptyset$}|S_j \subset C_j \in \mathcal{T}_n) = \frac{1}{n+1} \frac{Q([S_j])}{Q([W])}$$ and so on. The waiting times are independent of each other.\\
In general, one gets
\begin{equation}\label{eq: Lebenszeit-Mecke-Grundformel-ohne-Cowan}
\mathbb{P}\left(X_{S_j} \leq t|S_j \subset C_j \in \mathcal{T}_{n-1}=\mathcal{T}^s\right) = \sum_{k=n}^\infty \mathbb{P}\left(\sum_{i=n}^k T_i^M \leq t\right) \frac{1}{k} \frac{Q([S_j])}{Q([W])} \prod_{i=n}^{k-1} \left(1-\frac{1}{i} \frac{Q([S_j])}{Q([W])}\right).\end{equation}
Let us first calculate what one gets for $\mathbb{P}\left(\sum_{i=n}^k T_i^M \leq t\right)$ or the density of this respectively:
\begin{Lemma}
\label{Lemma: n-bis-k-te-Lebenszeit-kleiner-t}
The equation
\begin{equation}\label{eq: n-bis-k-te-Lebenszeit-kleiner-t}\mathbb{P}\left(\sum_{i=n}^k T_i^M \leq t\right) = \frac{1}{(k-n)!} \frac{k!}{(n-1)!} Q([W]) \int_0^t \left(e^{Q([W])x}-1\right)^{k-n} e^{-kQ([W])x} dx 
\end{equation} holds.
\end{Lemma}
\textbf{Proof}\\
We use the abbreviation $S_n^k = \sum_{i=n}^k T_i^M$. It is sufficient to show that for the density $f_{S_n^k}(x)$ of the probability distribution \begin{equation}\label{eq: Dichte-der-Summe-ab-n}f_{S_n^k}(x) = \frac{1}{(k-n)!} \frac{k!}{(n-1)!} Q([W]) \left(e^{Q([W])x}-1\right)^{k-n} e^{-kQ([W])x}\end{equation} holds.\\
The proof is by induction over $k$.\\
For the base case $k=n$, because of $T_n^M \sim \mathcal{E}(nQ([W]))$ the equation $f_{S_n^n} (x)= n Q([W]) e^{-nQ([W])x}$ should hold. Indeed,
$$f_{S_n^n}(x) = \frac{1}{(n-n)!} \frac{n!}{(n-1)!} Q([W]) \left(e^{Q([W])x}-1\right)^{n-n} e^{-nQ([W])x} = n Q([W]) e^{-nQ([W])x}.$$
Let now equation (\ref{eq: Dichte-der-Summe-ab-n}) be true for any $k$. Then for $k+1$, due to the convolution formula (the waiting times are independent of each other)
$$\begin{array}{rl}& f_{S_n^{k+1}}(x)\\&\\ = & \int_0^x \frac{1}{(k-n)!} \frac{k!}{(n-1)!} Q([W]) \left(e^{Q([W])u}-1\right)^{k-n} e^{-kQ([W])u}(k+1)Q([W]) e^{-(k+1)Q([W])(x-u)}du\\&\\
= & \frac{1}{(k-n)!} \frac{(k+1)!}{(n-1)!} (Q([W]))^2 \int_0^x \left(e^{Q([W])u}-1\right)^{k-n} e^{-kQ([W])u} e^{-(k+1)Q([W])x+kQ([W])u+Q([W])u}du\\&\\
= & \frac{1}{(k-n)!} \frac{(k+1)!}{(n-1)!} e^{-(k+1)Q([W])x} (Q([W]))^2\int_0^x \left(e^{Q([W])u}-1\right)^{k-n}  e^{Q([W])u}du\\&\\
= & \frac{1}{(k+1-n)!} \frac{(k+1)!}{(n-1)!} e^{-(k+1)Q([W])x} \frac{1}{Q([W])} (Q([W]))^2\int_0^x (k+1-n) Q([W]) \left(e^{Q([W])u}-1\right)^{k-n}  e^{Q([W])u}du\\&\\
= & \frac{1}{(k+1-n)!} \frac{(k+1)!}{(n-1)!} e^{-(k+1)Q([W])x} Q([W]) \left[\left(e^{Q([W])u}-1\right)^{k+1-n} \right]_{u=0}^{u=x} \\&\\
= & \frac{1}{(k+1-n)!} \frac{(k+1)!}{(n-1)!} e^{-(k+1)Q([W])x} Q([W]) \left(e^{Q([W])x}-1\right)^{k+1-n}
\end{array}$$
holds what is exactly what equation (\ref{eq: Dichte-der-Summe-ab-n}) yields for $k+1$. \hfill $\Box$\\
\begin{Lemma}
Let a time $s$ be fixed. At this time $s$, let us have a tessellation $\mathcal{T}^s$ with an arbitrary number of cells. Let a convex set $S_j$ be contained in the cell $C_j$ ($S_j \subset C_j \subset W$). For the waiting time $X_{S_j}$ for this convex set $S_j$ to be hit from the time $s$ on, \begin{equation}\label{eq: lifetime-convex-set}\mathbb{P}\left(X_{S_j} \leq t|S_j \subset C_j \in \mathcal{T}^s\right) = 1 - e^{-t Q([S_j])}\end{equation} holds.
\end{Lemma}
\textbf{Proof}\\
Let us first keep the number $n$ of quasi-cells fixed.
For equation (\ref{eq: Lebenszeit-Mecke-Grundformel-ohne-Cowan}) we get (at some point we will abbreviate $A = 1-e^{-Q([W])x}$) 
$$\begin{array}{rl}& \mathbb{P}\left(X_{S_j} \leq t|S_j \subset C_j \in \mathcal{T}_{n-1}=\mathcal{T}^s\right)\\&\\ 
= &\sum_{k=n}^\infty \mathbb{P}\left(\sum_{i=n}^k T_i^M \leq t\right) \frac{1}{k} \frac{Q([S_j])}{Q([W])} \prod_{i=n}^{k-1} \left(1-\frac{1}{i} \frac{Q([S_j])}{Q([W])}\right)\\&\\
= &\sum_{k=n}^\infty  \frac{1}{(k-n)!} \frac{k!}{(n-1)!} Q([W])\int_0^t \left(e^{Q([W])x}-1\right)^{k-n} e^{-kQ([W])x}dx \frac{1}{k} \frac{Q([S_j])}{Q([W])} \left[\prod_{i=n}^{k-1} \frac{1}{i} \left(i - \frac{Q([S_j])}{Q([W])}\right)\right]\\&\\
= &\sum_{k=n}^\infty  \frac{1}{(k-n)!} \frac{k!}{(n-1)!} \int_0^t \left(e^{Q([W])x}-1\right)^{k-n} e^{-kQ([W])x}dx Q([S_j]) \left( \prod_{i=n}^k \frac{1}{i} \right) \left[\prod_{i=n}^{k-1} \left(i - \frac{Q([S_j])}{Q([W])}\right)\right]\\&\\
= &\sum_{k=n}^\infty  Q([S_j]) \frac{1}{(k-n)!} \frac{k!}{(n-1)!} \int_0^t \left(e^{Q([W])x}-1\right)^{k-n} e^{-kQ([W])x}dx  \frac{(n-1)!}{k!} \left[\prod_{i=n}^{k-1} \left(i - \frac{Q([S_j])}{Q([W])}\right)\right]\\&\\
= &Q([S_j]) \int_0^t \sum_{k=n}^\infty   \frac{1}{(k-n)!}  \left(e^{Q([W])x}-1\right)^{k-n} e^{-kQ([W])x} \frac{\Gamma\left(k-\frac{Q([S_j])}{Q([W])}\right)}{\Gamma\left(n-\frac{Q([S_j])}{Q([W])}\right)} dx\\&\\
= &Q([S_j]) \int_0^t \sum_{k=0}^\infty   \frac{1}{k!}  \left(e^{Q([W])x}-1\right)^k e^{-(n+k)Q([W])x} \frac{\Gamma\left(n+k-\frac{Q([S_j])}{Q([W])}\right)}{\Gamma\left(n-\frac{Q([S_j])}{Q([W])}\right)} dx\\&\\\end{array}$$
$$\begin{array}{rl}
= &\frac{Q([S_j])}{\Gamma\left(n-\frac{Q([S_j])}{Q([W])}\right)} \int_0^t e^{-nQ([W])x} \sum_{k=0}^\infty   \frac{1}{k!}  \left(e^{Q([W])x}-1\right)^k \left(e^{-Q([W])x}\right)^k \Gamma\left(n+k-\frac{Q([S_j])}{Q([W])}\right) dx\\&\\
= &\frac{Q([S_j])}{\Gamma\left(n-\frac{Q([S_j])}{Q([W])}\right)} \int_0^t e^{-nQ([W])x} \sum_{k=0}^\infty   \frac{1}{k!}  \left(1-e^{-Q([W])x}\right)^k  \Gamma\left(n+k-\frac{Q([S_j])}{Q([W])}\right) dx\\&\\
= &\frac{Q([S_j])}{\Gamma\left(n-\frac{Q([S_j])}{Q([W])}\right)} \int_0^t e^{-nQ([W])x} \sum_{k=0}^\infty   \frac{1}{k!}  A^k  \int_0^\infty u^{n+k-\frac{Q([S_j])}{Q([W])}-1} e^{-u} du dx\\&\\
= &\frac{Q([S_j])}{\Gamma\left(n-\frac{Q([S_j])}{Q([W])}\right)} \int_0^t e^{-nQ([W])x} \int_0^\infty e^{-u} u^{n-\frac{Q([S_j])}{Q([W])}-1} \left(\sum_{k=0}^\infty   \frac{1}{k!}  A^k   u^k \right) du dx\\&\\
= &\frac{Q([S_j])}{\Gamma\left(n-\frac{Q([S_j])}{Q([W])}\right)} \int_0^t e^{-nQ([W])x} \int_0^\infty e^{-u} u^{n-\frac{Q([S_j])}{Q([W])}-1} e^{uA} du dx\\&\\
= &\frac{Q([S_j])}{\Gamma\left(n-\frac{Q([S_j])}{Q([W])}\right)} \int_0^t e^{-nQ([W])x} \int_0^\infty e^{-u(1-A)} u^{n-\frac{Q([S_j])}{Q([W])}-1}  du dx\\&\\
\stackrel{(a)}{=} &\frac{Q([S_j])}{\Gamma\left(n-\frac{Q([S_j])}{Q([W])}\right)} \int_0^t e^{-nQ([W])x} \int_0^\infty e^{-v} \left(\frac{v}{1-A}\right)^{n-\frac{Q([S_j])}{Q([W])}-1}  \frac{1}{1-A} dv dx\\&\\
= &\frac{Q([S_j])}{\Gamma\left(n-\frac{Q([S_j])}{Q([W])}\right)} \int_0^t e^{-nQ([W])x} (1-A)^{\frac{Q([S_j])}{Q([W])}-n} \int_0^\infty e^{-v} v^{n-\frac{Q([S_j])}{Q([W])}-1}  dv dx\\&\\
= &\frac{Q([S_j])}{\Gamma\left(n-\frac{Q([S_j])}{Q([W])}\right)} \int_0^t e^{-nQ([W])x} (1-A)^{\frac{Q([S_j])}{Q([W])}-n} \Gamma\left(n-\frac{Q([S_j])}{Q([W])}\right)dx\\&\\
= &Q([S_j]) \int_0^t e^{-nQ([W])x} (e^{-Q([W])x})^{\frac{Q([S_j])}{Q([W])}-n} dx\\&\\
= &Q([S_j]) \int_0^t e^{-nQ([W])x} e^{-Q([W])x\frac{Q([S_j])}{Q([W])}} e^{Q([W])xn} dx\\&\\
= &Q([S_j]) \int_0^t  e^{-xQ([S_j])}  dx\\&\\
= &Q([S_j]) \left[- \frac{1}{Q([S_j])} e^{-xQ([S_j])}\right]_{x=0}^{x=t}\\&\\
= & 1 - e^{-t Q([S_j])}.
\end{array}$$
Equation (a) follows from the substitution $v = u(1-A)$. \\
Let us now have a time $s$. The probability that there are exactly $n$ quasi-cells in the tessellation $\mathcal{T}^s$ (or that $\mathcal{T}^s = \mathcal{T}_{n-1}$) is just $$\mathbb{P}(\mathcal{T}^s = \mathcal{T}_{n-1}) = e^{-Q([W])t} \left(1 - e^{-Q([W])t}\right)^{n-1}.$$
Thus, we get 
$$\begin{array}{rl}\mathbb{P}\left(X_{S_j} \leq t|S_j \subset C_j \in \mathcal{T}^s\right) = & \sum_{n=1}^\infty \mathbb{P}(\mathcal{T}^s = \mathcal{T}_{n-1}) \left(1 - e^{-tQ([S_j])}\right)\\&\\
= & \sum_{n=1}^\infty e^{-Q([W])t} \left(1 - e^{-Q([W])t}\right)^{n-1} \left(1 - e^{-tQ([S_j])}\right)\\&\\ = & \left(1 - e^{-tQ([S_j])}\right) \sum_{n=1}^\infty e^{-Q([W])t} \left(1 - e^{-Q([W])t}\right)^{n-1}\\&\\
= & 1 - e^{-tQ([S_j])}.\end{array}$$ Thus the lemma is proven. \hfill $\Box$\\
It is worth to mention that (as shown by the last equation) the result does \textit{not} depend on the number of quasi-cells $n$ at the time $s$. For homogeneous $Q$, this result is the same result one has for the STIT process. Obviously, the lifetime of a cell $C_j$ (which is a convex set contained within a fixed cell, namely $C_j$) can be described in this manner as well.

\section{Proofs of Mecke's Conjectures}
\label{sec: Proofs-Conjectures-1-2}
\subsection{Conjecture 1}
Lemma \ref{Lemma: X-t-in-equally-likely} makes clear the relation between those properties Mecke calls 'characteristics' and the waiting time in a state with $n$ quasi-cells. Let us have a tessellation in a window $\hat{W}$ with characteristics $\hat{Q}$ and $\hat{t}$; then we get $$\mathbb{P}(\hat{N}_t = k) = e^{-\hat{t} } \left(1-e^{-\hat{t} }\right)^k.$$

\begin{Lemma} (Mecke's Conjecture 1)
\label{Lemma: Beweis-Conjecture-1}
Let $\mathcal{T}_W$ be a mixed line-generated tessellation in $W$ with characteristics $Q$ and $tQ([W])$, and let $\hat{W}$ be a window with $\hat{W} \subset W$ and $Q([\hat{W}]) > 0$. Then the cutout of $\mathcal{T}_W$ in $\hat{W}$ can be interpreted as a mixed line-generated tessellation in $\hat{W}$ with characteristics $$\hat{Q} = \frac{1}{Q([\hat{W}])} Q(\cdot \cap [\hat{W}]) \textrm{   and   } \hat{t} = tQ([\hat{W}]).$$
\end{Lemma}
\textbf{Proof}\\
We will first examine the tessellation $\mathcal{T}_W$ in $W$ with characteristics $Q$ and $tQ([W])$. For the probability that until a time $t$ the convex subset $\hat{W} \subset W$ was hit, according to equation (\ref{eq: lifetime-convex-set}) $$\mathbb{P}(X_{\hat{W}} \leq t) = 1 - e^{-Q([\hat{W}])t}$$ holds.\\
If we condition on $\hat{W}$ being hit the hitting line has distribution $\frac{1}{Q([\hat{W}])} Q(\cdot \cap [\hat{W}]) = \hat{Q}$.\\
Let us now examine the tessellation $\mathcal{T}_{\hat{W}}$ in $\hat{W}$ with characteristics $\hat{Q}=\frac{1}{Q([\hat{W}])} Q(\cdot \cap [\hat{W}])$ and $\hat{t} = tQ([\hat{W}])$. The distribution of the number of decisions $\hat{\nu}(t)$ until time $t$ is $$\mathbb{P}(\hat{\nu}(t) = k) = e^{-\hat{t}} \left(1 - e^{-\hat{t}}\right)^k = e^{-t Q([\hat{W}])} \left(1-e^{-t Q([\hat{W}])}\right)^k.$$ 
From this, we can deduce the lifetime of the first cell $\hat{W}$ to be $$\mathbb{P}(X_{\hat{W}} \leq t) = 1-e^{-t Q([\hat{W}])}.$$
Thus, the distribution of the lifetimes of $\hat{W}$ is the same in $\mathcal{T}_W$ and $\mathcal{T}_{\hat{W}}$; additionally, the distribution of the segment dividing $\hat{W}$ is identical as well.\\
Let us now have $\mathcal{T}_W \cap \hat{W} = \mathcal{T}_{\hat{W}}$ at an arbitrary time. Then, in $\mathcal{T}_W$ there exist the cells  $Cells(\mathcal{T}_W) = \{C^W_1, ..., C^W_n\}$ and accordingly in $\mathcal{T}_{\hat{W}}$ the cells $Cells(\mathcal{T}_{\hat{W}}) = \{C^W_1 \cap \hat{W}, ..., C^W_n \cap \hat{W}\} \setminus \{\emptyset\}$. Note that some of the intersections $C^W_j \cap \hat{W}$ can be empty; therefore the empty set is taken out of the set in order to have only 'real' cells with non-empty interior in $Cells(\mathcal{T}_{\hat{W}})$.\\
We now examine a cell $C \in Cells(\mathcal{T}_W)$ with $C \cap \hat{W} \neq \emptyset$. This cell has, as calculated above, the lifetime $X_C \sim \mathcal{E}(Q([C]))$. If we take a look at this cell's intersection with the subwindow $\hat{W}$ we have a waiting time $X_{C \cap \hat{W}} \sim \mathcal{E}(Q([C \cap \hat{W}]))$ for this convex set to be hit. For the distribution of the dividing line we have, due to the conditioning on the division of the set, $\frac{1}{Q([\hat{W}])} Q([C \cap \hat{W}]) \frac{1}{Q([C \cap \hat{W}])} Q([\cdot \cap C \cap \hat{W}]) = \hat{Q}(\cdot \cap [C]).$\\
In the tessellation $\mathcal{T}_{\hat{W}}$, we have a cell $\hat{C} = C \cap \hat{W}$ with a lifetime $X_{\hat{C}} \sim \mathcal{E}(Q([\hat{C}])) = \mathcal{E}(Q([C \cap \hat{W}]))$. For the distribution of the line dividing $\hat{C}$ we have $\hat{Q}(\cdot \cap [\hat{C}])$.\\
So, the waiting time for the set $C \cap \hat{W}$ in $\mathcal{T}_W$ to be hit and the lifetime of the cell $C \cap \hat{W}$ in $\mathcal{T}_{\hat{W}}$ respectively are identically distributed. Because of $\hat{Q}(\cdot \cap [C]) = \hat{Q}(\cdot \cap [C \cap \hat{W}]) = \hat{Q}(\cdot \cap [\hat{C}])$ the distributions of the dividing lines are identical as well.\\
Under the condition of the equality $\mathcal{T}_W \cap \hat{W} = \mathcal{T}_{\hat{W}}$ the distributions of the time of the next segment falling in $\hat{W}$ are identical in both considered windows as are the distributions of that next segment. As we always start in the same configuration of an empty subset $\hat{W} \subset W$ and window $\hat{W}$ respectively, the lemma is proven.\hfill $\Box$


\subsection{Conjecture 2}
With Lemma \ref{Lemma: Beweis-Conjecture-1}, Mecke's Conjecture 2 can be proven in quite a straightforward manner:
\begin{Lemma} (Mecke's Conjecture 2)
\label{Lemma: Beweis-Conjecture-2} 
The class of all mixed line-generated tessellations (related to $Q$) as a whole is stable under iteration in the following sense: Every operation of iteration maps the mentioned class into itself, i.e. an iterated mixed line-generated tessellation is again a mixed line-generated tessellation. If the mixed line-generated tessellation $\mathcal{T}^t$ is iterated according to the law $P^s$ of $\mathcal{T}^s$, then the law $P^t \boxplus P^s$ of the outcome fulfils $$P^t \boxplus P^s = P^{t+s}.$$
\end{Lemma}
\textbf{Proof}\\
Let $\mathcal{T}^t$ be a tessellation with the cells $Cells(\mathcal{T}^t)=\{Z_1, ..., Z_\kappa\}$. Then each of those cells $Z_j$ has a lifetime $X_{Z_j} \sim \mathcal{E}(Q([Z_j]))$ independent from the other lifetimes which because of the memorylessness of the exponential distribution does non depend on the time the cell was created before the time $t$. After the lifetime has expired (provided it is smaller than $s$), a segment of a line with distribution $\frac{1}{Q([Z_j])} Q(\cdot \cap [Z_j])$ falls into the cell. Afterwards, the process goes on with its cells and their exponentially-distributed lifetimes until time $s$. Thus, one gets the resulting tessellation $\mathcal{T}^{t+s}$.\\
According to Lemma \ref{Lemma: Beweis-Conjecture-1}, one can interpret the cutout $\mathcal{T}_W \cap Z_j$ of $\mathcal{T}_W$ with characteristics $Q$ and $sQ([W])$ as a process $\mathcal{T}_{Z_j}$ with characteristics $\hat{Q} = \frac{1}{Q([Z_j])} Q(\cdot \cap [Z_j])$ and $\hat{s} = s Q([Z_j])$. If one considers a cell $Z_j$ now, its lifetime is $\mathcal{E}(Q([Z_j]))$-distributed as verified in the proof of Lemma \ref{Lemma: Beweis-Conjecture-1}; after this lifetime's expiry, a segment falls with the corresponding line having a distribution $\hat{Q} = \frac{1}{Q([Z_j])} Q(\cdot \cap [Z_j])$.\\
This cutout process runs independently in all cells $Z_j, j=1, ..., \kappa,$ with the same lifetime and segment distribution as in the process $\mathcal{T}^{t+s}$. Thus, because the processes are identically distributed, $$P^t \boxplus P^s = P^{t+s}$$ holds as claimed in Mecke's Conjecture 2.\hfill $\Box$

\section{Conclusion}
With these two results from Lemma \ref{Lemma: Beweis-Conjecture-1} and Lemma \ref{Lemma: Beweis-Conjecture-2} and the result from \cite[Theorem 1]{Conjecture-3}, all conjectures by Mecke stated in \cite{Mecke-inhomogen} have been proven.

\begin{center}\textbf{Acknowledgement}\end{center}
I am thankful to Joseph Mecke for providing me with a link to the paper \cite{Faltung-Summe-exponentialverteilter-Zufallsgroessen}. I am once again very thankful to Werner Nagel for his help in putting together this paper.

\bibliographystyle{plain} \bibliography{literatur}

\end{document}